\documentclass[10pt,reqno]{article}
\usepackage{latexsym, amsmath, amssymb, a4, epsfig}
\usepackage{graphicx}
\usepackage{algorithm,algpseudocode}
\usepackage{multirow, array}

\usepackage[normalem]{ulem}
\usepackage{epstopdf}
\usepackage{xcolor}
\newcommand\ruline{\bgroup\markoverwith{\textcolor{red}{\rule[-0.7ex]{2pt}{0.8pt}}}\ULon}

\newcommand{\nm}{\noalign{\smallskip}}
\newcommand{\ds}{\displaystyle}

\newcommand{\p}{\partial}

\newcommand{\eqnref}[1]{(\ref {#1})}

\newcommand{\Rbb}{\mathbb{R}}

\newcommand{\la}{\langle}
\newcommand{\ra}{\rangle}

\newcommand{\Kcal}{\mathcal{K}}

\newcommand{\Dcal}{\mathcal{D}}

\newcommand{\Gd}{\delta}

\newcommand{\Gvf}{\varphi}

\newcommand{\Gs}{\sigma}

\newcommand{\GD}{\Delta}

\newcommand{\GG}{\Gamma}

\newcommand{\GO}{\Omega}

\newcommand{\beq}{\begin{equation}}
\newcommand{\eeq}{\end{equation}}
\newcommand{\beqa}{\begin{eqnarray}}
\newcommand{\eeqa}{\end{eqnarray}}

\numberwithin{equation}{section}
\numberwithin{figure}{section}

\begin{document}
\title{A non-iterative method for the electrical impedance tomography based on joint sparse recovery\thanks{\footnotesize This work is supported by the Korean Ministry of Education, Sciences and Technology through NRF grants Nos. 2010-0017532 (to H.K), NRF-2013R1A1A3012931 (to M.L), NRF-2009-0081089 (to O.K.L. and J.C.Y.), the R\&D Convergence Program of NST (National Research Council of Science \& Technology) of Republic of Korea (Grant CAP-13-3-KERI) (to O.K.L. and J.C.Y.).}}

\author{Ok Kyun Lee\thanks{Department of Bio and Brain Engineering, Korea Advanced Institute of Science and Technology, Daejeon 305-701, S. Korea (okkyun2@kaist.ac.kr, jong.ye@kaist.ac.kr).} \and Hyeonbae Kang\thanks{Department of Mathematics, Inha University, Incheon 402-751, S. Korea (hbkang@inha.ac.kr).} \and Jong Chul Ye\footnotemark[2] \and Mikyoung Lim\thanks{Department of Mathematical Sciences, Korea Advanced Institute of Science and Technology, Daejeon 305-701, S. Korea (mklim@kaist.ac.kr).}}
\maketitle

\begin{abstract}
The purpose of this paper is to propose a non-iterative method for the inverse conductivity problem of recovering multiple small anomalies from the boundary measurements.
 When small anomalies are buried in a conducting object, the electric potential values inside the object can be expressed by integrals of densities with a common sparse support on the location of anomalies. Based on this integral expression, we formulate the reconstruction problem of small anomalies as a joint sparse recovery and present an efficient non-iterative recovery algorithm of small anomalies. Furthermore, we also provide a slightly modified algorithm to reconstruct an extended anomaly.
We validate the effectiveness of the proposed algorithm over the linearized method and the MUSIC algorithm by numerical simulations.
 \end{abstract}



\section{Introduction}

Electrical impedance tomography (EIT) is a noninvasive imaging technique to reconstruct the electrical property of a medium based on the boundary measurement of the voltages that result from the injected currents. The electrical properties of a material are characterized by
\beq\label{eq:ElecProp}
\tilde{\sigma} = \sigma + j\omega \epsilon,
\eeq
where $\sigma$ is the electric conductivity, $\omega$ the angular frequency of the applied current waveform, and $\epsilon$ the electric permittivity. With the help of relatively low-cost imaging equipment and the fact that various materials such as biological tissues, certain rocks and fluids have their own $\sigma$ and $\epsilon$ values \cite{Borcea2002EIT}, EIT has been applied for various clinical and industrial applications such as monitoring internal organs in human body, finding mineral deposits on earth, and nondestructive inspection \cite{Alessandrini1998EIT,Costa2009EIT,Parker1984EIT}.
However, the inverse problem of EIT is nonlinear and ill-posed due to the nonlinear coupling of the electrical potential to the electrical material properties and due to the compactness of the forward mapping.

To circumvent the nonlinearity, one-step linearization methods or iterative methods are commonly used in practice \cite{Cheney1999EIT,Saulnier2001EIT}. A one-step linearization method employs the computed value of the internal electrical potential corresponding to the background electrical property, instead, in the place of the true, but unknown, internal potential. While this procedure is quite fast, it produces a non-negligible error in the reconstruction. On the other hand, an iterative method gives more accurate result by updating the solution multiple times, but it now suffers from
the ill-posedness nature in the EIT problem besides the computational burden of solving the forward problem in each updating step.
To overcome the ill-posedness in the inverse problem of EIT, it has been studied extensively by a variety of techniques. Especially for the problem of electrical anomalies detection, various algorithms have been proposed, among which are the small volume expansion method \cite{AK01, expan,  FV89}, the projection algorithm \cite{KSY02}, the simple pole algorithm \cite{KL04}, the linear sampling method \cite{BHV}, the MUSIC (MUltiple SIgnal Classification) algorithm \cite{BHV}, and the topological derivative algorithm \cite{book_new}. It is worth mentioning that the MUSIC algorithm was proposed by Schmidt \cite{Sc86} for source separation in signal theory (see also \cite{devaney2}) and has been applied to many imaging problems for which we refer to \cite{book_new}. The topological derivative algorithm was proposed by Eschenauer et al \cite{esch} for shape optimization and has been successfully applied to imaging problems of diverse contexts including the anomalies detection. Even though various reconstruction methods have been proposed, most of them are either linear approximations or iterative methods.

In this paper, we propose a non-iterative reconstruction method for EIT problem which can be exact if the measurement is noiseless. The algorithm exploits the joint sparsity of the induced current source on the anomaly for different current injection directions. The joint-sparse recovery method comes from the compressed sensing theory which deals with the under-determined linear problem to recover sparse signals that share common non-zero support \cite{chen2006trs,Kim2010CMUSIC}.

In the EIT problem, the idea of sparsity has been employed for the regularization in the optimization based inversion. In \cite{KKSV,MI}, the regularization for least squares problems with the $l^1$ prior was considered. Recently, the sparsity regularization has been proposed to reconstruct the conductivity distribution when the object under consideration has a sparse representation with respect to a certain basis \cite{GKSM,JZM}. However, these sparsity based regularizations are also within the framework of the conventional reconstruction methods, {\it i.e.}, either linearized or iterative methods.

One of the important contributions of the proposed method is that, by exploiting the joint sparsity, one can obtain an accurate reconstruction of the anomalies without linearization or iteration. The idea to exploit the joint sparsity in anomaly detection problem originally comes from the previous researches that overcome the nonlinearity in the inverse problem of diffuse optical tomography \cite{Lee2011CDOT,Lee2013DOT}, and now we apply it under the circumstance of EIT problem. More specifically, we change the non-linear EIT inverse problem to the joint sparse recovery problem and obtain the electrical anomalies by following the three simple steps. First, the common non-zero support of induced currents due to the presence of anomalies and the values of induced currents are obtained by the joint sparse recovery method. Secondly, the internal electric potential is estimated. Finally, the electrical properties are calculated by solving the associated linear problem.
Numerical results of the proposed method for both sparse and extended anomalies are provided and are compared to those of a linear approximation and the MUSIC algorithm in order to validate the efficiency of the proposed method.

The paper is organized as follows. Mathematical background of EIT is given in Section \ref{sec:Background}, and a brief introduction to compressed sensing and joint sparse recovery is given in Section \ref{sec:CS}. The proposed method using joint sparsity is described in detail in Section \ref{sec:Proposed}, which is followed by its implementation in Section \ref{sec:Implementation}. Section \ref{sec:Results} provides numerical results, and finally, a discussion and conclusion is provided in Section \ref{sec:Conclusion}.

\section{Problem Formulation}
\label{sec:Background}

In this paper, we deal with the static problem ($\omega=0$) for simplicity.
Let $\GO$ be a bounded domain in $\Rbb^d$, $d=2,3$, which is occupied by the homogeneous material of the conductivity 1. We suppose that a finite number of isotropic anomalies are embedded in the background domain $\GO$. We denote the anomalies by $D_1, \ldots, D_N$ and the corresponding conductivities by $\Gs_1,\ldots,\Gs_N$ where $N$ is the number of anomalies. So the conductivity profile of the domain $\GO$ is given by
\beq
\ds\Gs= \ds\chi(\GO \setminus \cup_{p=1}^N D_p) + \sum_{p=1}^N \Gs_p \chi(D_p),
\eeq
where the symbol $\chi(D)$ indicates the characteristic function of $D$.
We let $\Gs_p$\rq{}s be smooth and satisfy the ellipticity condition
$$0<c\leq \Gs_p(x)\leq C<\infty\quad \mbox{for } x\in D_p, \ \ p=1,\dots,N$$
for some positive constants $c$ and $C$.

The EIT problem we consider in this paper is to find the locations, geometric features, and conductivities of anomalies using a finite number of pairs of voltages (Dirichlet data) and currents (Neumann data) on the boundary of $\GO$. Let $g_1, \ldots, g_M \in L^2_0(\p\GO)$ be the $M$ number of  given currents on $\partial \GO$. Here $L^2_0(\partial \GO)$ means the set of square integrable functions defined on $\p\GO$ with zero means. The currents $g_k$'s are given functions and will be considered as column vectors in the associated linear system that will be explained later. The corresponding internal potential $u_k$ in $\GO$ for $1\leq k\leq M$ satisfies the Neumann boundary value problem
\beq\label{eqn:u_k}
 \left \{
 \begin{array}{ll}
 \ds \nabla  \cdot ( \Gs(x)  \nabla u_k) =0 \quad \mbox{in } \GO,  \\
 \nm \ds \frac{\p u_k}{\p\nu} \Bigr |_{\p\GO} =g_k,\\
 \nm \ds \int_{\p\GO} u_k \, d\Gs=0.
 \end{array}
 \right .
\eeq
We tackle the reconstruction problem of the anomalies $D_p$'s and their conductivities $\sigma_p$'s based on an integral representation formula of $u_k$. We derive the formula in the remaining of this section.

The Neumann function $N(\cdot,y)$ on $\GO$ is the solution to
\beq\label{Neumann}
 \left \{
 \begin{array}{ll}
 \ds - \Delta_x N(x, y) = \Gd_y(x) \quad &\mbox{in } \GO, \\ \nm \ds
 \frac{\p }{\p\nu_x} N(x,y)\Big|_{\p\GO} = - \frac{1}{|\partial \GO|} , \\
 \nm \ds
 \int_{\p\GO} N(x,y) \, d\sigma(x) =0
 \end{array}
 \right .
\eeq
for $y\in \GO$.
Let $U_k$ be the electric potential in absence of anomalies, {\it i.e.}, the solution to
\beq\label{eqn:U}
 \left \{
 \begin{array}{l}
 \ds \Delta U_k =0 \quad \mbox{in } \GO,  \\
 \nm \ds \frac{\p U_k}{\p\nu} \Bigr |_{\p\GO} =g_k, \\
 \nm \ds \int_{\p\GO} U_k \, d\Gs=0.
 \end{array}
 \right .
\eeq
Then $U_k$ can be represented as
$$
U_k(x) = \int_{\p\GO} N(x,y) g_k(y) ~d\Gs(y),\quad x \in \GO.
$$
Note that because of the third condition in \eqnref{eqn:u_k} and the second equation in \eqnref{Neumann} we have
$$
\int_{\p\GO} \frac{\p}{\p\nu_y} N(x,y) u_k(y)~ d\Gs(y)=0,\quad x \in \GO.
$$
Thus we have
$$
U_k(x) = \int_{\p\GO}  \left[ N(x,y) \frac{\p u_k}{\p\nu}(y) - \frac{\p}{\p\nu_y} N(x,y) u_k(y) \right] \;d\Gs(y).
$$
We then have from the Green's identity
\begin{align*}
U_k(x) & = - \int_{\GO\setminus \cup_{p=1}^N D_p} \GD_y N(x,y) u_k(y)~ dy \\
& \qquad + \sum_{p=1}^N\int_{\p D_p}  \left[ N(x,y) \frac{\p u_k}{\p\nu}\Big|_+ (y) - \frac{\p}{\p\nu_y} N(x,y) u_k(y) \right] \;d\Gs(y) .
\end{align*}
It then follows from the transmission conditions (continuity of flux and potential) of $u_k$ along $\p D_p$'s that
\begin{align*}
U_k(x) & = - \int_{\GO\setminus \cup_{p=1}^N D_p} \GD_y N(x,y) u_k(y)~ dy \\
& \quad + \sum_{p=1}^N \int_{\p D_p}  \left[\Gs_p(y)  N(x,y) \frac{\p u_k}{\p\nu}\Big|_- (y) - \frac{\p}{\p\nu_y} N(x,y) u_k(y) \right] \;d\Gs(y) \\
& = - \int_{\GO\setminus \cup_{p=1}^N D_p} \GD_y N(x,y) u_k(y)~ dy - \int_{\cup_{p=1}^N D_p} \GD_y N(x,y) u_k(y)~ dy \\
&\quad +\int_{\cup_{p=1}^N D_p}  N(x,y) \GD u_k(y)~ dy+ \sum_{p=1}^N \int_{\p D_p}(\Gs_p(y)-1)  N(x,y) \frac{\p u_k}{\p\nu}\Big|_- (y) ~d\Gs(y) \\
& =u_k(x) +  \int_{\cup_{p=1}^N D_p} (\Gs(y)-1)\nabla_y N(x,y) \cdot \nabla u_k(y) ~dy.
\end{align*}
Here the symbols $-$ and $+$ indicate the limits from inside and outside of $D_p$ to $\p D_p$, respectively.
Finally we obtain the following formula:
\beq\label{ukxUkx}
u_k(x)-U_k(x) = \int_{\cup_{p=1}^N D_p}   (1-\Gs(y))\nabla_y N(x,y) \cdot \nabla u_k(y)~ dy, \quad x \in \overline{\GO}.
\eeq

The joint sparsity method of this paper is based on the formula \eqnref{ukxUkx}. However, it is not easy to compute the Neumann function unless domains are disks or balls. So, we apply the Calder\'on preconditioner \cite{book2}.
For that let $\GG(x)$ be the fundamental solution to the
Laplacian, {\it i.e.},
$$
\GG (x) =
\begin{cases}
\ds \frac{1}{2\pi} \ln |x|\;, \quad & d=2 \;, \\ \nm \ds
-\frac{1}{4\pi} |x|^{-1}\;, \quad & d = 3 \;,
\end{cases}
$$
and define an operator (called the Neumann-Poincar\'e operator) by
\beq\notag
\Kcal_{\p\GO} [\Gvf] (x) =
\frac{1}{\omega_d} \int_{\partial \GO} \frac{\la
 y -x, \nu_y \ra}{|x-y|^d} \Gvf(y)\,d\sigma(y)\;, \quad x \in \p \GO.
\eeq
Then it is known (see \cite{book2}) that
\beq\notag
\left( -\frac{1}{2} I + \Kcal_{\p\GO} \right ) \big[N(\cdot, y)\big](x)= \GG(x-y) \ \mbox{modulo constant}, \ x \in \p\GO, \ y \in \GO.
\eeq
By applying the operator $-\frac{1}{2} I + \Kcal_{\p\GO}$ to both sides of the equality in \eqnref{ukxUkx}, we obtain
\beq\label{ukxUkx2}
\left( -\frac{1}{2} I + \Kcal_{\p\GO} \right ) \big[(u_k - U_k)|_{\p\GO}\big](x) = \int_{\cup_{p=1}^N D_p} (1-\Gs(y)) \nabla_y \GG(x-y) \cdot \nabla u_k(y)\; dy, \quad x \in \p\GO.
\eeq
So the problem is to reconstruct $D_p$ (and $\Gs_p$) using $\left( -\frac{1}{2} I + \Kcal_{\p\GO} \right ) [(u_k - U_k)|_{\p\GO}]$ for $k=1, \ldots M$.

We emphasize that \eqnref{ukxUkx2} holds for $x\in \p\GO$. In order to derive an integral relation which holds for $x \in \GO$, we may use the double layer potential which is defined to be
\beq\label{doublelayer}
\Dcal_{\p\GO} [\Gvf] (x) := \int_{\p \GO} \frac{\p}{\p\nu_y} \Gamma (x-y) \Gvf (y) \, d\sigma(y) \;, \quad x \in \Rbb^d \setminus \p\GO .
\eeq
Then because of the jump relation
\beq\notag
\Dcal_{\p\GO} [\Gvf] \big |_- (x) = \biggl( \frac{1}{2} I + \Kcal_{\p\GO} \biggr) [\Gvf] (x),
\quad x \in \p\GO,
\eeq
we obtain
\beq\label{eq:ukxUkxInt}
u_k(x)-U_k(x)=\Dcal_{\p\GO} \big[(u_k-U_k)|_{\p\GO}\big](x)+\int_{\cup_{p=1}^N D_p} (\sigma(y)-1)\nabla_y \Gamma(x-y)\cdot\nabla u_k(y)\;dy,\quad x\in\GO.
\eeq

In linear approximation approaches, one may use a further approximation of the formula \eqnref{ukxUkx} using smallness of $D_p$:
The unknown potential $u_k$ is the small perturbation of $U_k$ and \eqnref{ukxUkx} becomes
\beq\notag
u_k(x)-U_k(x) \approx\int_{\cup_{p=1}^N D_p}   (1-\Gs(y))\nabla_y N(x,y) \cdot \nabla U_k(y)~ dy, \quad x \in \overline{\GO}.
\eeq
By taking $-\frac{1}{2} I + \Kcal_{\p\GO}$ to the boundary value of $u_k-U_k$, we have
\beq\label{eqn:born}
\left( -\frac{1}{2} I + \Kcal_{\p\GO} \right ) \big[(u_k - U_k)|_{\p\GO}\big](x) \approx \int_{\cup_{p=1}^N D_p}(1-\Gs(y)) \nabla_y \GG(x-y) \cdot \nabla U_k(y)\; dy, \quad x \in \p\GO.
\eeq
And if $z_p$ represents the location of $D_p$, then we have a simpler approximation formula
\beq\label{ukxUkx3}
u_k(x) - U_k(x) \approx \sum_{p=1}^N \nabla_y N(x, z_p) \cdot M_p \nabla U_k(z_p), \quad x \in \p\GO.
\eeq
Here $M_p$ is a $d \times d$ matrix associated with $D_p$ and $\sigma_p$, and is called the polarization tensor (PT). The formula \eqnref{ukxUkx3} was first found in \cite{FV89} when the conductivity of the inclusions is $\infty$. See \cite{CY03,book2} for derivation of the formula.

We will implement in the following sections the joint sparse recovery method using the formula \eqnref{ukxUkx2} and \eqnref{eq:ukxUkxInt} in order to recover the anomalies and their conductivities. Then we will compare the performance with two conventional methods, one based on \eqnref{eqn:born} (the linearized EIT method) and one based on \eqnref{ukxUkx3} (the MUSIC algorithm).
It is worth making a comment on \eqnref{ukxUkx2} and \eqnref{eq:ukxUkxInt}. In order to use these formula we need to have $(u_k-U_k)|_{\p\GO}$, namely, the measurement $u_k$ and $U_k$ on the whole boundary $\p\GO$. The background domain $\Omega$ is assumed to be known so that the solution $U_k$ in absence of anomalies can be computed.  If measurements on only a part of the boundary are available, these formula cannot be used, and instead one may use \eqnref{ukxUkx}. Additionally, we will discuss the anomalies reconstruction based on joint sparse recovery for the smooth background conductivity and the electrode model of EIT in Section 7.

\section{Compressed Sensing}
\label{sec:CS}

This section briefly introduces the compressed sensing theory and addresses the joint sparse recovery problem to make a seamless flow from the previous section to the next. Compressed sensing theory is the state of the art in the field of signal processing that enables the recovery of the signal beyond the Nyquist limit based on the sparsity of the signal \cite{CaRoTa06}. As an example, let us consider the under-determined linear system of $y=Ax$ that has many solutions. When the signal $x$ has a sparsity, the accurate recovery of the signal is possible using the compressed sensing theory as described in the following problem \cite{CaRoTa06}:
\begin{equation}\label{eq:y=Ax0}
\min_{x} \|x\|_0, \mbox{ subject to } y = Ax,
\end{equation}
where $y \in \mathbb{R}^{m\times 1}$, $A \in \mathbb{R}^{m \times n}$, and $x \in \mathbb{R}^{n\times 1}$ with $m<n$. Here, $\|x\|_0$ denotes the number of non-zero elements in the vector $x$, and $x$ having a sparsity means that $||x||_0$ is much smaller than $n$. The uniqueness of the solution to the problem \eqref{eq:y=Ax0} is guaranteed by the following condition \cite{DoEl03}:
\begin{equation}\label{eq:y=Ax0uniq}
\|x\|_0 < \frac{\mbox{spark}(A)}{2},
\end{equation}
where $\mbox{spark}(A)$ is the smallest possible number $\ell$ such that there exist $\ell$ linearly dependent columns of $A$. However, \eqref{eq:y=Ax0} is an $NP$-hard problem that every possible combination of supports should be considered. Therefore, the following $l_1$ minimization is widely used in practice \cite{Candes2008IntroCS}:
\begin{equation}\label{eq:y=Ax1}
\min_{x}\|x \|_1, \mbox{ subject to } y = Ax,
\end{equation}
where $\|\cdot\|_1$ denotes the $l_1$ norm. The original $NP$-hard problem is now relaxed into a convex optimization problem. The beauty of compressed sensing is that \eqref{eq:y=Ax1} provides the exactly same solution as \eqref{eq:y=Ax0} if the so called restricted isometry property (RIP) is satisfied \cite{Candes2008IntroCS}. It has been shown that for many class of random matrices, the RIP is satisfied with extremely high probability if the number of measurement satisfies  $m\geq c k\log (n/k)$, where $k=\|x\|_0$ and $c$ is a positive constant \cite{Candes2008IntroCS}.

The multiple measurement vector (MMV) problem \cite{chen2006trs,Kim2010CMUSIC} is a generalization of the single measurement vector (SMV) problem defined in \eqref{eq:y=Ax0}. The MMV problem is the signal recovery problem  to exploit a set of sparse signal vectors that share common non-zero supports, in other words, a set of signal vectors that have a joint sparsity. Specifically, let $||X||_0$ denote the number of rows that have non-zero elements in the matrix $X$. Then, the MMV problem addresses the following:
\begin{equation} \label{eq:y=Ax0MMV}
\min_{X} ||X||_0, \quad \mbox{ subject to } Y = AX,
\end{equation}
where $Y \in \mathbb{R}^{m \times M}$, $X \in \mathbb{R}^{n \times M}$, and $M$ denotes the number of measurement vectors. Intuitively, we can tell at a glance that the MMV problem \eqref{eq:y=Ax0MMV} contains more information than the SMV problem \eqref{eq:y=Ax0}, so that it provides better reconstruction results unless the column vectors in $X$ are all about the same. Theoretically, \eqref{eq:y=Ax0MMV} has the unique solution if and only if
\begin{equation}\label{eq:y=Ax0MMVuniq}
||X||_0 < \frac{\mbox{spark}(A)+ \mbox{rank}(Y)-1}{2} \left(\leq \frac{m+\mbox{rank}(Y)}{2}\right),
\end{equation}
where $\mbox{rank}(Y)$ denotes the rank of $Y$ and it may increase with the number of measurement vectors \cite{chen2006trs,Davies2012MMV}. Note that $\mbox{rank}(Y)$ term in \eqref{eq:y=Ax0MMVuniq} clearly shows the advantage of the MMV problem over the SMV problem. {Since calculating $\mbox{spark}(A)$ is also an $NP$-hard problem, it is hard to know the real least upper bound of the recoverable number of non-zero elements. We can roughly estimate it using the number of measurement points as in the rightmost inequality in \eqref{eq:y=Ax0MMVuniq} instead. According to the definition of $\mbox{spark}(A)$, it is obvious that the equality is satisfied when all the combination of $m$ columns of $A$ are independent.}

\begin{figure}[htbp]
\centering\includegraphics[width=13cm]{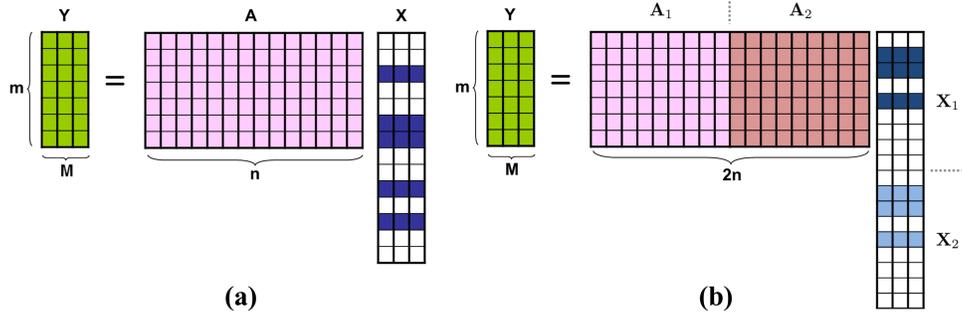}
\caption{Joint sparsity model in (a) general and (b) EIT problem in two dimensions.}
\label{fig:CS}
\end{figure}

There are various kinds of joint sparse recovery algorithms to solve the MMV problem including the convex relaxation \cite{cotter2005ssl, Kim2010CMUSIC, tropp2006ass2,wipf2007ebs}. Fig.\;\ref{fig:CS}(a) illustrates the general joint sparse recovery problem, but the EIT problem we consider can be addressed with a bit special pairwise joint sparse recovery problem (Fig.\;\ref{fig:CS}(b)) as will be described in the next section.

\section{Recovery Method using Joint Sparsity}
\label{sec:Proposed}
This section describes the non-iterative exact reconstruction method for EIT problem in $\Rbb^2$ based on the joint sparsity recovery algorithm. It reconstructs the anomalies accurately without an iterative procedure or a linear approximation. Even though we only deal with two dimensions (in this and following sections), it is straightforward to extend it to three dimensions.

Let us first fix some notations. For the potential function $u_k$ for $1\leq k\leq M$, we define the current on the anomalies as
\beq\label{eq:Ijk}
I_{k}(y)=\begin{cases}
(1-\sigma_p(y))\nabla u_k(y)\quad&\mbox{for } y\in D_p,\ p=1,\dots,N,\\
0\quad&\mbox{for }y\in\GO\setminus\cup_{p=1}^N D_p.
\end{cases}
\eeq
Then the formula \eqnref{ukxUkx2} becomes
\beq\label{eq:ukxUkx4}
\left( -\frac{1}{2} I + \Kcal_{\p\GO} \right ) \big[(u_k - U_k)|_{\p\GO}\big](x) = \int_\GO \nabla_y \GG(x-y) \cdot I_{k}(y) \;dy, \quad x \in \p\GO.
\eeq
Remind that the anomalies $D_1, \ldots, D_p$ are located at fixed positions despite of the different boundary conditions $g_1, \ldots, g_M$, whereas the currents $I_k$ on the anomalies vary. Therefore, assuming the sparsity for the support set $\cup_{p=1}^N D_p$, the problem \eqnref{eq:ukxUkx4} is a joint sparsity problem since the non-zero current location (the non-zero rows in the associated linear equation) is independent of the boundary currents (the given columns).
To describe it more specifically, let us assume that $I_{k}$, $1\leq k\leq M$, is approximated by either piecewise constant functions or splines: \beqa\label{eq:CurrentModel}
\ds I_{k}(y) = \left[
        \begin{array}{c}
           \ds \sum\limits_{j=1}^{n} I_{k,1}(y_{(j)} )\hskip .5mm b_1 ( y,y_{(j)} ) \\
          \ds  \sum\limits_{j=1}^{n} I_{k,2} ( y_{(j)} )\hskip .5mm b_2 ( y,y_{(j)} )
        \end{array}
            \right], \quad y\in\GO,
\eeqa
where $b_d$, $d=1,2$, is the basis function for the $d$-th coordinate and  $\{ y_{(j)} \}_{j=1}^{n}$ is the finite discretization points of $\GO$.

After substituting \eqref{eq:CurrentModel} into \eqref{eq:ukxUkx4}, we have the following equation:
\beq\label{eq:ukxUkx5}
\ds\left( -\frac{1}{2} I + \Kcal_{\p\GO} \right ) \big[(u_k - U_k)|_{\p\GO}\big](x)
= \sum_{d=1}^2 \sum_{j=1}^n \tilde{\GG}_d(x,y_{(j)}) \hskip .5mm I_{k,d}(y_{(j)}),\quad x \in \p\GO,
\eeq
where
$$
\tilde{\GG}_d(x,y_{(j)}) =  \int_{\GO} \nabla_{y,d} \GG(x-y)\hskip .5mm b_d( y,y_{(j)} )\; dy.$$
Here $\nabla_{y,d} \GG(x-y)$ ($d=1,2$) means the $d$-th coordinate component of $\nabla_{y} \GG(x-y)$.
We now define, respectively, the sensing matrix $A=[A_1,A_2]\in \mathbb{R}^{m\times 2n}$, the currents $X=[X_1^T,X_2^T]^T\in \mathbb{R}^{2n\times M}$ and the measurements $Y\in \mathbb{R}^{m\times M}$
as \begin{gather*}
\ds\big(A_d\big)_{ij}=\tilde{\GG}_d(x_{(i)},y_{(j)}),\quad\big(X_d\big)_{jk}=I_{k,d}(y_{(j)}),\quad\mbox{and } \big(Y\big)_{ik}=\left( -\frac{1}{2} I + \Kcal_{\p\GO} \right )\big[(u_k - U_k)|_{\p\GO}\big](x_{(i)}),\end{gather*}
where $\{ x_{(i)} \}_{i=1}^{m}$ is the collection of the finite number of measurement locations on $\p\GO$. Then we can formulate \eqref{eq:ukxUkx5} as the following matrix equation:
\beq\label{eq:MatEq}
Y = AX+E = [A_1, A_2]\left[
        \begin{array}{c}
            X_1 \\
            X_2
        \end{array}
            \right] + E
\eeq
with the measurement noise $E \in \mathbb{R}^{m\times M}$. From \eqnref{eq:Ijk}, the solution $X$ to \eqref{eq:MatEq} has a pairwise joint sparsity meaning that $X_1$ and $X_2$ are nonzero at the same rows which correspond to the positions where the anomalies are located, as is illustrated in Fig.\;\ref{fig:CS}(b).
Based on this equation, we can formulate the following joint sparse recovery problem \cite{chen2006trs}:
\beq\label{eq:JSM}
\min_{X}~||X||_0,~~\mbox{subject to}~~ ||Y - AX||_F^2 \leq \epsilon,
\eeq
where $||\cdot||_F$ denotes the Frobenius norm.

There are well known algorithms to solve the joint sparse recovery problem. After having the solution $X$ to \eqnref{eq:JSM} by applying one of those algorithms, one can estimate the anomaly positions by collecting $y_{(j)}$ whose corresponding currents $(X_d)_{jk}$ are nonzero for all $k$ and $d$. Since this criterion reconstructs the points which belong to any one of $D_p$'s, we denote the obtained anomaly positions by $\hat{D}$. With the solution $X$ to \eqnref{eq:JSM} and $\hat{D}$, we have the current on the anomaly as well, say $\hat{I}_{k}(y_{(j)})$ for $y_{(j)}\in\hat{D}$. Then the unknown solution $u_k$ can be now easily estimated using \eqref{eq:ukxUkxInt} as
\beq\label{eq:ukxUkxInt2}
\hat{u}_k(x) = U_k(x) + \Dcal_\GO \big[(u_k-U_k)|_{\p\GO}\big](x) - \int_{\hat{D}}\nabla_y \Gamma(x-y)\cdot \hat{I}_{k}(y)\; dy,\quad x\in \hat{D}.
\eeq
Finally, the conductivity $\sigma$ is calculated by solving
\beq\label{eq:ukxUkx6}
\left( -\frac{1}{2} I + \Kcal_{\p\GO} \right ) \big[(u_k - U_k)|_{\p\GO}\big](x) =- \int_{\hat{D}}(\sigma(y) -1) \nabla_y \GG(x-y) \cdot \nabla \hat{u}_k(y)\; dy, \quad x \in \p\GO
\eeq
for $1\leq k\leq M$. Note that every term in \eqnref{eq:ukxUkx6} except $(\sigma -1)$ is now known since the potential $u_k$, which was initially measured only on $\p\GO$, is estimated on the whole anomalies $\hat{D}$ from \eqnref{eq:ukxUkxInt2}. We emphasize that \eqref{eq:ukxUkx6} is a linear equation for $(\sigma(y_{(j)})-1)$.  Hence neither linear approximation nor the iterative update is required. Furthermore, we can expect more efficient and less ill-posed reconstruction procedure due to the knowledge of the estimated position of anomalies.

\section{Implementations}
\label{sec:Implementation}

\subsection{Joint Sparse Recovery}\label{subsec:JointSparseRecovery}

To solve the problem \eqref{eq:JSM}, we use the multiple sparse Bayesian learning (M-SBL) algorithm \cite{wipf2007ebs}. The M-SBL algorithm assumes that the noise element and $X$ follow the {\it i.i.d.} Normal distribution with ${\rm vec}(E) \sim \mathcal{N}(0,\lambda I_{mM})$ and ${\rm vec}(X)\sim \mathcal{N}(0, I_M \otimes \Gamma)$, respectively, where $I_M$ denotes the $M\times M$ identity matrix and $\Gamma$ is the common variance component for the $i$-th row values of $X$ given by the diagonal matrix of entries $\gamma_i$'s. Considering the pairwise joint sparsity in \eqref{eq:JSM}, we further assume $\gamma_i=\gamma_{i+n}$ for $i=1,2,\cdots,n$. It is worth mentioning that even though the original M-SBL is derived from the Bayesian framework, recent theoretical analysis \mbox{\cite{Wipf2011Latent}} shows that M-SBL is indeed a sparse recovery algorithm that can be used in deterministic framework without assuming any statistics of $X$. More specifically, it solves the following problem:
\beq\label{eq:MSBLcost}
\min_{X} \|Y-AX\|_F^2 + \lambda f(X),
\eeq
where the penalty function $f(X)$ is given as
\beq\label{eq:MSBLcost2}
\min_{\gamma_i\geq 0} {\rm Tr}(X'\Gamma^{-1}X) + M \log \det (A\Gamma A' + \lambda I_m).
\eeq
{Rather than the constraint optimization problem in \eqref{eq:JSM}, the unconstrained form of the cost function in \eqref{eq:MSBLcost} is usually used to deal with noisy measurement.  Here, the regularization parameter $\lambda$ is determined based on the noise level. Furthermore, the penalty term of M-SBL in \eqref{eq:MSBLcost2} is shown to impose the sparsity more effectively compared to the conventional $l^p$ norm as shown in \mbox{\cite{Wipf2011Latent}}.}

The step-by-step procedure for M-SBL is summarized in Algorithm \ref{alg:PseudocodeMSBL}. The lines from $7$ to $9$ in the algorithm is the pruning step for the variance components inherent in the M-SBL, and $M$ in the lines $6$ and $10$ is the number of the measurement vectors on $\p \GO$.
We normalized the sensing matrix $A$ for each column to have a unit $l^2$ norm before applying Algorithm \ref{alg:PseudocodeMSBL}.

\begin{algorithm}
\caption{Pseudocode implementation of the M-SBL.}
\label{alg:PseudocodeMSBL}
\begin{algorithmic}[1]
\State Set $\mbox{Iter}_{\max}$.
\State Set $k=0$ and $\lambda^{(0)}=0.01 \times \sigma_{\max}(A)^2$.
\Statex Set $\gamma^{(0)}_i=1$ for $i=1,2,\cdots,2n$ and $\Gamma^{(0)}=\mbox{diag}(\mathbf{\gamma}^{(0)})$.
\For{$k=1,\dots, \mbox{Iter}_{\max}$}
    \State Set $\Lambda=\left(A\Gamma^{(k-1)}A' + \lambda^{(k-1)}I \right)^{-1}$
    \State Update $X^{(k)}=\Gamma^{(k-1)}A'\Lambda Y$.
    \State Update $\gamma_i^{(k)}=\gamma_{i+n}^{(k)}=\sqrt{\frac{||X^{i(k)}||^2+||X^{i+n(k)}||^2}{M\left(A_i'\Lambda A_i + A_{i+n}'\Lambda A_{i+n}\right)}}$ for $i=1,2,\cdots,n$ and set $\Gamma^{(k)}$.
    \If{${\gamma_i^{(k)}}/{{\rm max}(\mathbf{\gamma}^{(k)})} < 10^{-3}$}
        \State $\gamma_i^{(k)}=0$
    \EndIf
    \State Update $\lambda^{(k)} = \sqrt{\frac{||Y-AX^{(k)}||_F^2}{M\cdot\mbox{Tr}\left(\Lambda \right)}}$.
\EndFor
\end{algorithmic}
\end{algorithm}

As a result of the M-SBL algorithm, we obtain the solution $X$ to \eqnref{eq:JSM}. We now obtain the support of $(1-\sigma)$ by collecting $y_{(j)}$ whose corresponding currents $I_{k,d}(y_{(j)})$ are nonzero for all $k$ and $d$. More precisely, we estimate the position of anomalies based on the following criterion:
\beq\label{eq:SpecMSBL2}
\hat{D} = \left\{y_{(j)} \;\bigg\vert\; \frac{p_{(j)}}{\max \;(p)} > \epsilon \right\},
\eeq
where $\epsilon$ is a threshold and $p$ the spectrum of the current defined by
\beq\label{eq:SpecMSBL}
p_{(j)} = \sqrt{\sum_{d=1}^{2} \sum_{k=1}^{M} \left[\left(X_d\right)_{jk}\right]^2}\quad\mbox{for }1\leq j\leq n.
\eeq
The current's spectrum and relative errors with various thresholding values $\epsilon$ for sparse and extended target examples are provided in Fig.\;\ref{fig:ResultEpsilon}.


One can apply a preconditioning procedure before the M-SBL algorithm to deal with the ill-posedness in the inversion of the sensing matrix which is inherited from the ill-posedness nature in the inverse problem of EIT.
Let us denote the singular value decomposition of the sensing matrix $A$ as $A=USV'$. Then the regularized preconditioning matrix $P$ becomes $P = (S^2+ \lambda I)^{-1/2}U'$ and the preconditioned problem for \eqref{eq:MatEq} can be restated as
\beq\label{eq:MatEqPre}
PY = PAX + PE.
\eeq
See \cite{Jin2012precond}. We now apply the M-SBL algorithm with $PY$ and $PA$ instead of $Y$ and $A$.

\subsection{Conductivity Recovery}\label{subsec:Born}

After solving the joint sparse recovery problem in \eqref{eq:JSM} and finding the unknown value of $u_k(x)$ using \eqref{eq:ukxUkxInt2}, we can calculate the conductivity of the anomalies from the linear equation for $(1-\sigma)$ in \eqref{eq:ukxUkx6}. Let us denote by $\{\hat{y}_{(j)}\}_{j=1}^{\tilde{n}}$ the estimated points of $\hat{D}$ from \eqref{eq:SpecMSBL2} and $\delta$ the area of the discretized grid, then the discretized version of \eqref{eq:ukxUkx6} is as follows:
\beq\label{eq:CondRecon}
y=    \left[
        \begin{array}{c}
            y_1 \\  \vdots  \\ y_M
        \end{array}
    \right] =
           \left[
        \begin{array}{c}
            A_1 \\  \vdots \\ A_M
        \end{array}
    \right]
   \left[
        \begin{array}{c}
            x_1 \\ \vdots \\ x_{\tilde{n}}
        \end{array}
    \right] = Ax,
\eeq
where $y_k\in \mathbb{R}^{m\times 1}$, $A_k\in\mathbb{R}^{m\times\tilde{n}}$, and $x\in \mathbb{R}^{\tilde{n}\times 1}$ are given by
$$\big(y_{k}\big)_i = \left( -\frac{1}{2} I + \Kcal_{\p\GO} \right ) \big[(u_k - U_k)|_{\p\GO}\big](x_{(i)}),\quad\left(A_k\right)_{i,j}=-\nabla_y \GG(x_{(i)}-\hat{y}_{(j)}) \cdot \nabla \hat{u}_k (\hat{y}_{(j)})\; \delta,$$ and $x_{j}=\sigma(\hat{y}_{(j)})-1$, respectively. Note that the number of unknowns in the discretized domain is reduced from $n$ to $\tilde{n}$, and the sensing matrix $A$ is accurate if the estimates of $\hat{y}_{(j)}$'s and $\hat{u}_k$'s are precise. To solve \eqnref{eq:CondRecon}, we use the following constrained optimization problem with $l_1$ penalty for noise robust reconstruction:
\beqa\label{eq:CondRecon2}
 \mbox{arg}\min \limits_x & ||x||_{1}\\
 \mbox{subject to} & || Ax-y ||_2 \leq \epsilon. \nonumber
\eeqa

\begin{figure}[htbp]
\centering\includegraphics[width=9cm,height=3cm]{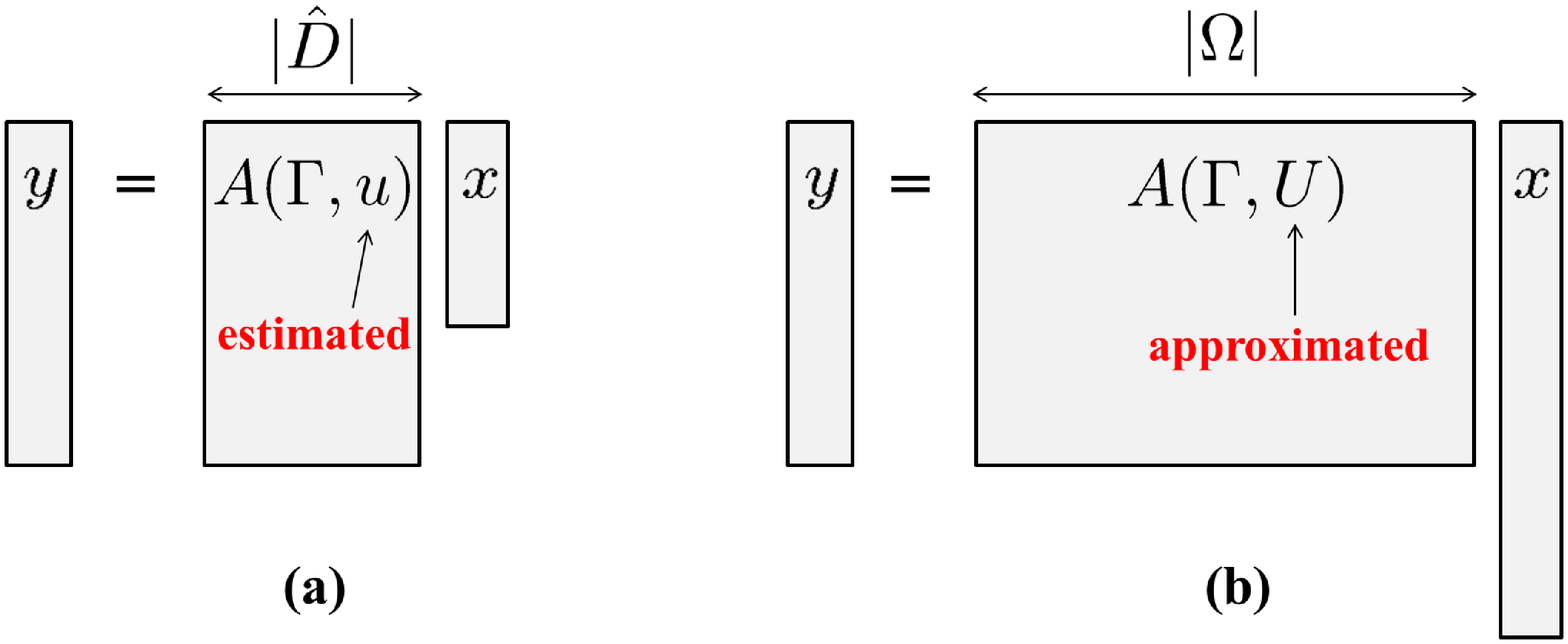}
\caption{Conductivity recovery problem of the (a) proposed method and (b) linearized EIT problem.}
\label{fig:Compare}
\end{figure}

We can construct a linear system similar to \eqnref{eq:CondRecon} based on the linearized EIT problem \eqnref{eqn:born}.
In other words, the sensing matrix $A$ in \eqref{eq:CondRecon} becomes $A_k\in\mathbb{R}^{m\times n}$, $1\leq k\leq M$, such that
$$\big(A_k\big)_{ij}=-\nabla_y \GG(x_{(i)}-y_{(j)}) \cdot \nabla U_k (y_{(j)})\; \delta\quad\mbox{with }y_{(j)}\in \GO.$$
The size of sensing matrix $A$ is increased back from $\tilde{n}$ to $n$ because of lack of information on the location of anomalies and the accuracy of $A$ becomes worse because of the approximation error. The advantage of the linearized EIT problem over the conventional iterative methods is in the computing speed. However, the proposed method using joint sparsity has even higher computational efficiency as illustrated in Fig.\;\ref{fig:Compare}(a) compared to the linearized EIT problem in Fig.\;\ref{fig:Compare}(b). Moreover, the proposed method has also better accuracy in the anomalies recovery. We compare the speed and the accuracy of the two methods for various examples in the next section.

 To solve the problem \eqref{eq:CondRecon2}, we exploit a constrained split augmented Lagrangian shrinkage algorithm (C-SALSA) \cite{Afonso2011CSALSA} whose pseudocode implementation is described in Algorithm \ref{alg:PseudocodeCSALSA}. We normalized the sensing matrix $A$ for each column to have a unit $l_2$ norm before applying the C-SALSA algorithm.
In the pseudocode, the constant $\tau$ is set to be $\tau=c\bar{s}_0^{(1)}$ at the $0$-th iteration, where $\bar{s}_0^{(1)}$ is the average of $|s_0^{(1)}|$ and $c$ is the constant chosen manually to be optimal among $c=\{8, 4, 2, 1, 1/2, 1/4, 1/8\}$. The parameter $\mu$ is the continuation factor, which is chosen to be $\mu=1.001$ for the {linearized EIT problem} and $\mu=1.01$ for the proposed method in this paper.
The $l_1$ minimization problem changes into the simple soft thresholding by the proximal mapping
$$\Psi_{\tau g_1}(s_k^{(1)}) = \mbox{arg}\min\limits_v \frac{1}{2}||v-s_k^{(1)}||_2^2 + \tau ||v||_1$$ with $v$  such as component-wise operation of $$v_j = \mbox{sign}\big([s_k^{(1)}]_j\big)\max\big\{\big|[s_k^{(1)}]_j\big|-\tau,\;0\big\}\quad\mbox{for } j\mbox{-th component}.$$ We made a criterion such that $\Psi_{\tau g_i}$ is the Moreau proximal mapping of $\tau g_i$ with $g_1(x) = ||x||_1$ and $g_2 = \iota_{E(\epsilon,y)}$, where $\iota_{E(\epsilon,y)}$ is the indicator function of the $\epsilon$-radius Euclidean ball centered at $y$.  Here $\epsilon$ is a constant given by $\epsilon=c||y||_2$ with $c$ manually chosen to be optimal among $c=\{0.02, 0.04, 0.06, 0.08, 0.1, 0.2, 0.3\}$.
As a stopping criterion, we use the relative change of the cost function in \eqref{eq:CondRecon2} and perform the algorithm until $\left| \left(C_{k}-C_{k-1}\right)/C_{k} \right| < 10^{-8}$ is satisfied, where $C_{k}$ is the cost function at the $k$-th iteration.

 \begin{algorithm}
\caption{Pseudocode implementation of C-SALSA for the EIT problem in \eqref{eq:CondRecon2}.}
\label{alg:PseudocodeCSALSA}
\begin{algorithmic}[1]
\State Set $k=0$, $H^{(1)}=I$, $H^{(2)}=A$, and choose $\tau>0$.
\State Set $v_0^{(i)}=d_0^{(i)}=0$, for $i=1,2$, except $v_0^{(2)}=y$.
\Repeat
\For{$i=1,2$}
    \State $\zeta_k^{(i)} = v_k^{(i)} + d_k^{(i)}$
\EndFor
\State $u_{k+1} = \left[ \sum\limits_{j=1}^{2}\left(H^{(j)}\right)^TH^{(j)} \right]^{-1} \sum\limits_{j=1}^{2}\left(H^{(j)}\right)^T\zeta_k^{(j)}$
\For{$i=1,2$}
    \State $v_{k+1}^{(i)} = \Psi_{\tau g_i} \left( s_k^{(i)} \right)$, where $s_k^{(i)} = H^{(i)} u_{k+1} - d_k^{(i)}$
    \State $d_{k+1}^{(i)} = d_{k}^{(i)} - H^{(i)} u_{k+1} + v_{k+1}^{(i)}$
\EndFor
\State $\tau = \tau / \mu$
\State $k = k + 1$
\Until{some stopping criterion is satisfied}
\end{algorithmic}
\end{algorithm}

{The sensing matrix of the proposed method and that of the linearized method are different, so the selected optimal parameters are distinct as will be shown later. Moreover, $\bar{s}_0^{(1)}$ values are also different since the explicit form of $s_0^{(1)}$, which is $\left(A^TA+I\right)^{-1}A^Ty$, depends on the sensing matrix.}


\subsection{Extended Target Recovery}

{The proposed method we described in the previous section consists of three steps. First, the target location and corresponding current values are reconstructed from the joint sparse recovery. Second, the unknown potential is estimated, and conductivities are calculated as a final step. While the proposed method is designed aiming to recover sparse anomalies, the other two steps are unrelated to the size of the target. Therefore, the proposed method can also be applied to the recovery of non-sparse targets as long as one can solve the (modified) first step. Unfortunately, it turns out that the estimation of the potential in the second step becomes incorrect for an extended target if we use the currents as well as the anomalies support obtained from the M-SBL in the first step. To alleviate this challenge for the non-sparse target, we suggest to solve the following linearized equation to estimate the currents distribution using only the estimated $\hat{D}$ from the M-SBL algorithm:}
\beq\label{eq:ukxUkx7}
\left( -\frac{1}{2} I + \Kcal_{\p\GO} \right ) \big[(u_k - U_k)|_{\p\GO}\big](x) = \int_{\hat{D}} \nabla_y \GG(x-y) \cdot I_k(y) \; dy, \quad x \in \p\GO.
\eeq
{We validate the proposed reconstruction method for extended target with the simulation results in the following section. }

\section{Numerical Simulations}
\label{sec:Results}

In this section we present numerical simulation results using the proposed method and compare them with those using the linearized method described in Section~\ref{subsec:Born} and also with the MUSIC algorithm.
For the sake of simplicity we abbreviate the linearized method explained in Section \ref{subsec:Born} by the linearized EIT.

\begin{figure}[htbp]
\centering\includegraphics[width=13cm]{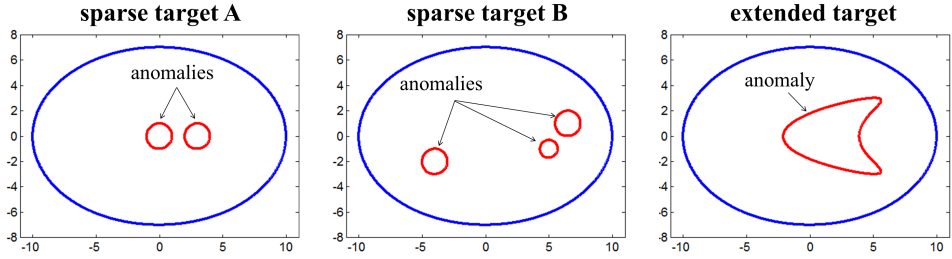}
\caption{Simulation geometry of sparse and extended targets.}
\label{fig:FOV}
\end{figure}

We show two sparse target examples and one extended target example. In all examples, $\GO$ is an ellipse of semi-major and semi-minor axes 10 and 7, and the background conductivity is $1$. The sparse anomalies in the first example (named the sparse target A) are the two unit disks away from $\p\GO$ whose distance is 1.  The sparse anomalies in the second example (named the sparse target B) are three disk shaped anomalies of different size placed arbitrarily. Lastly, we let the anomaly be a kite shaped extended target to test for the non-sparse target. The three examples are illustrated in Fig.\;\ref{fig:FOV}.

The left anomaly in the first example has the conductivity value of $2$ and the right one has the value of $5$. In the second example, conductivities of anomalies are $0.5$, $5$, and $2$ from left to right. In the third example,  the conductivity of the extended target is $5$. The field of view is discretized to have a grid size $0.5$ for reconstruction.

\begin{figure}[!htbp]
\centering\includegraphics[width=10cm]{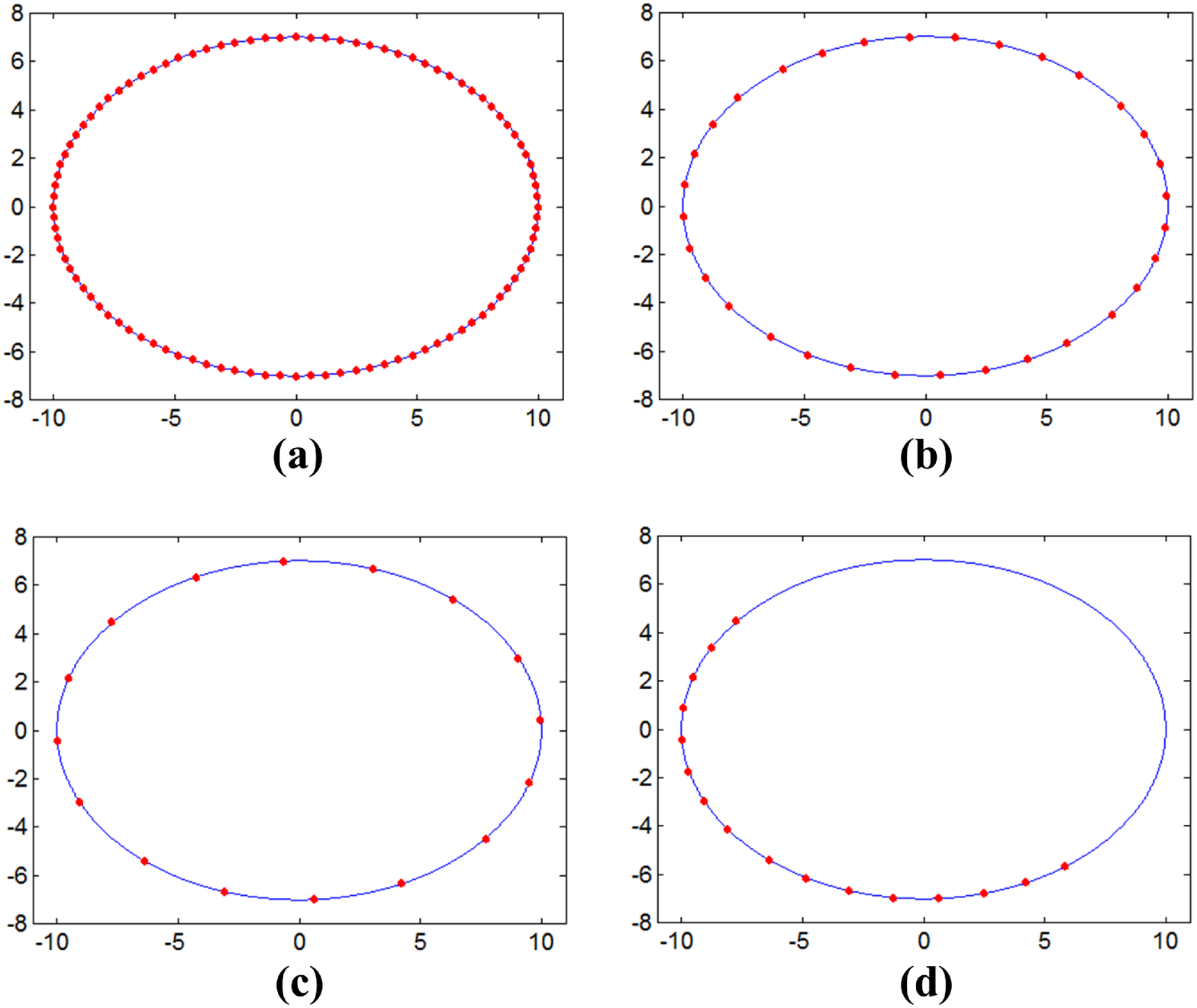}
\caption{Measurement points geometry on the boundary of $\GO$. Measurement points are $m$ number of uniformly sampled either on $\p \GO$ with (a) $m=100$, (b) $m=32$, (c) $m=16$, or on a half part of $\p \GO$ with (d) $m=16$. We denote the geometry (d) by $m=16p$ to indicate the partial measurements. }
\label{fig:FOVm}
\end{figure}

For the boundary current $g_k$ in \eqnref{eqn:u_k} we use \beq\label{eqn:g_k}g_k=\nabla H_k\cdot\nu_\GO\quad\mbox{for } 1\leq k\leq M,
\eeq
where $\nu_\GO$ is the outward unit normal vector to $\p\GO$ and
 $$H_1(x)=x_1,\ H_2(x)=x_2,\ H_3(x)=x_1^2-x_2^2,\ H_4(x)=x_1x_2$$ with $x=(x_1,x_2)$.
Either two measurements $(M=2)$ or four measurements $(M=4)$ are used in the numerical examples.
In order to acquire the measurement data $(u_k-U_k)\big|_{\p\GO}$, we solve the boundary integral equation \eqnref{eqn:u_k} and \eqnref{eqn:U} numerically. We use the computation method modified from that in \cite{book_new}, which is based on the expression of $u_k$ and $U_k$ in terms of the single layer potentials. See \cite{book_new} for the details of the numerical code.
It is worth to remark that the background solution $U_k$ for $g_k$ given by \eqnref{eqn:g_k} is actually $H_k$.
While $2000$ nodal points are used on each $\partial\Omega$ and $\partial D_j$ in the direct solver, the data only at a limited number of points on $\p\Omega$ are used in the reconstruction procedure where the sampling points number $m$ is either $100$, $32$, or $16$.
 See Fig.\;\ref{fig:FOVm} for the location of sampling points on $\p \GO$. In Fig.\;\ref{fig:FOVm}(d), $16$ measurement points are located on only part of the boundary and we let $m=16p$ refer to this case. For the reconstruction with $m=16p$, we follow the recovery process with \eqnref{ukxUkx} instead of \eqnref{eq:ukxUkxInt}. Hence it is necessary to compute the Neumann function $N(x,y)$ for all $y\in\GO$ for $m=16p$.
 Gaussian noise with a $\mbox{SNR}$ of $40dB$ was added to the boundary measurement vectors $(u_k-U_k)\big|_{\p\GO}$ for all simulations.

\subsection{Parameters}
M-SBL is used to estimate the position and current values of sparse anomalies; and for the extended anomaly, M-SBL is used only for the position and the truncated singular value decomposition (T-SVD) \cite{Hansen1987TSVD} is used for estimating current values based on \eqref{eq:ukxUkx7}. The regularization parameter $\lambda$ in the preconditioning in \eqref{eq:MatEqPre} for the extended target recovery is selected to be $\lambda=10^{-3}\times\sigma_{\max}(A)^2$, where $\sigma_{\max}(A)$ denotes the maximum singular value of $A$.

\begin{figure}[!htbp]
\centering\includegraphics[width=15cm]{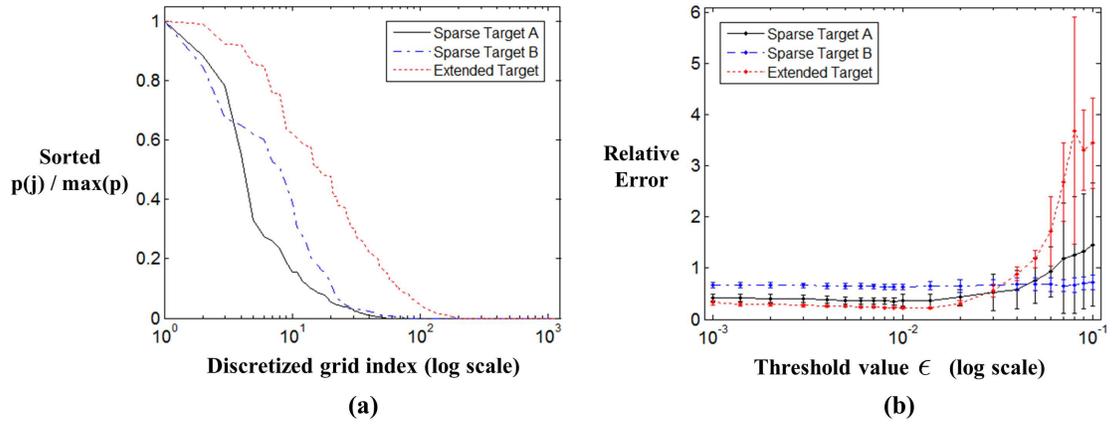}
\caption{(a) Normalized spectrum of the current in \eqref{eq:SpecMSBL} and (b) the relative error along the various thresholding values $\epsilon$ for the spectrum. The error bars indicate the standard deviation from $20$ trials when $m=100$ and $M=2$.}
\label{fig:ResultEpsilon}
\end{figure}

{As described in \eqref{eq:SpecMSBL2} and \eqref{eq:SpecMSBL}, the proposed method requires the threshold value $\epsilon$ for the recovery support of $(\sigma -1)$. Fig.\;\ref{fig:ResultEpsilon}(a) shows the normalized graph for the current's spectrum defined in \eqref{eq:SpecMSBL} for sparse and extended targets. We chose the threshold value $\epsilon=1\times 10^{-2}$ to cover the most of the dominant portion while avoiding the meaningless area at the same time. To make an observation of how the threshold value influences to the proposed method, we plot the error bars with varying the threshold values in Fig.\;\ref{fig:ResultEpsilon}(b), when $m=100$ and {$M=2$}.
For the quantitative measure, we calculate the relative error defined as
\beq\label{eq:relError}
\mbox{error}=\frac{||x_{true}-x_{recon}||^2}{||x_{true}||^2},
\eeq
where $x_{true}=\sigma_{true}-1$, $x_{recon}=\sigma_{recon}-1$, {and $||\cdot||^2$ is the squared $l_2$ norm}.
Fig.\;\ref{fig:ResultEpsilon}(b) shows that the under-estimation of the position of the anomalies deteriorates more severely than the over-estimation.}

The selected optimal parameters for reconstruction algorithms are summarized in Table\;\ref{tbl:ParameterSimul}.
These parameters are used for all examples except for the case of extended target with $m=16p$ where we use the increased value of $\epsilon=0.2||y||_2$ due to the deteriorated data fidelity term.
\begin{table}[!htbp]
\footnotesize
\vspace*{0.2cm}
\centerline{
\renewcommand{\arraystretch}{1.3}
\begin{tabular}{|c|c|c|c|c|} \hline
                & Algorithm & Sparse target A  & Sparse target B  & Extended target \\
\hline \hline
\multirow{4}{*}{Proposed Method}&\multirow{2}{*}{M-SBL}& $\mbox{Iter}_{\max}=15$   & $\mbox{Iter}_{\max}=14$   & $\mbox{Iter}_{\max}=9$ \\
                & & w/o preconditioning   & w/o preconditioning   & with preconditioning \\
\cline{2-5}
                &\multirow{2}{*}{C-SALSA} & $\tau= 1\bar{s}_0^{(1)}$  & $\tau= 0.5\bar{s}_0^{(1)}$  & $\tau=0.125 \bar{s}_0^{(1)}$ \\
                & & $\epsilon=0.04||y||_2$ & $\epsilon=0.02||y||_2$ & $\epsilon=0.06||y||_2$ \\
\hline
\multirow{2}{*}{Linearized EIT} & \multirow{2}{*}{C-SALSA}& $\tau=8 \bar{s}_0^{(1)}$ & $\tau=0.125 \bar{s}_0^{(1)}$ & $\tau=0.25 \bar{s}_0^{(1)}$\\
 & & $\epsilon=0.06||y||_2$ & $\epsilon=0.08||y||_2$ & $\epsilon=0.1||y||_2$\\
\noalign{\hrule height 0.5pt}
\end{tabular}
}\caption {Parameters used in simulation.}
\label{tbl:ParameterSimul}
\normalsize
\end{table}

\subsection{MUSIC algorithm}
{For the MUSIC algorithm, we apply the eigenvalue-based MUSIC algorithm as suggested in \cite{Mosher1992MEGmusic}. More specifically, MUSIC spectrum for $y_{(j)}\in\Omega$ is calculated based on the formulation in \eqref{ukxUkx3}:
\beq\label{eq:MUSIC}
\frac{1}{\lambda_{\min} \left( U_j' P U_j \right)},
\eeq
{where $\lambda_{\min}(X)$ denotes the minimum eigenvalue of $X$, $P$ is the orthogonal projector on the noise subspace of the measurement vectors, and $U_j$ denotes the two most-dominant eigenvectors of $\left[ (N_1)_j, (N_2)_j \right] \left[ (N_1)_j, (N_2)_j \right]'$. Here, $(N_1)_j$ and $(N_2)_j$ are the $j$-th columns of $N_1$ and $N_2$, respectively, where $N_k$ is composed of as following:}
$$\big(N_k\big)_{ij}= \nabla_{y,k} N(x_{(i)},y_{(j)}), \quad\mbox{with }x_{(i)}\in \p\GO,~y_{(j)}\in \GO$$
for $k=1,2$.

\subsection{Simulation results}

Before showing the conductivity reconstructions, we visualize the internal potentials that are used for the conductivity recovery. Remind that while $u_k$'s are estimated by the joint sparse recovery in the proposed method, the background potential $U_k$'s are used in the linearized EIT. We show $U_k$ and $u_k$ for $k=2$ near $\hat{D}$ in Fig.\;\ref{fig:ResultU} at the 2nd and 3rd columns, respectively. For both sparse and extended targets, the internal potentials obtained from the proposed method are similar to the true
$u_k$ (1st column).
\begin{figure}[!htbp]
\centering\includegraphics[width=14cm]{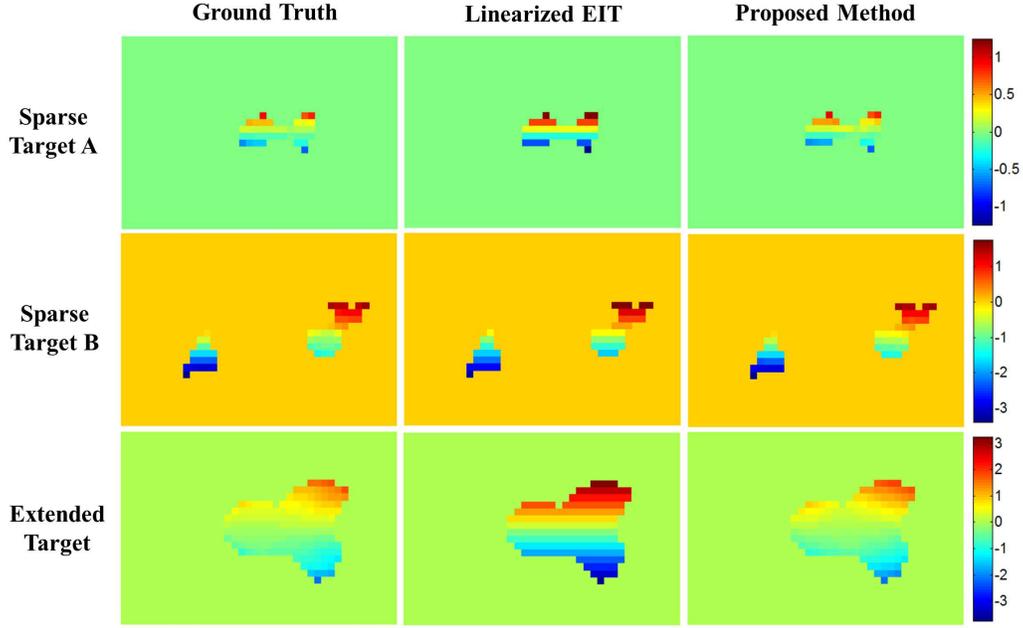}
\caption{True $u_2(y)$ (1st column), $U_2(y)$ for the linearized EIT (2nd column), and estimated $u_2(y)$ using the proposed method (3rd column) for $y \in \hat{D}$.}
\label{fig:ResultU}
\end{figure}

In Fig.\;\ref{fig:ResultAll}, we illustrate reconstructions of $(\sigma-1)$ for sparse and extended targets from the proposed and the linearized EIT method using two measurements $g_1$ and $g_2$, {\it i.e.}, $M=2$. We show examples with the measurement point geometry $m=100$ explained in Fig.\;\ref{fig:FOVm} for the linearized EIT and with various geometries $m=100,32,16,16p$ for the proposed method.
When $m=100$, two anomalies in the sparse target A are clearly separated from each other in the result by the proposed method and have distinct conductivity values which are close to the ground truth better than those by the linearized EIT. As the case of sparse target B is more complicated than that of A, the overall reconstruction performance is downgraded. However, the proposed method still provides comparable results to ground truth except the underestimated conductivity value of the centered {anomaly} while right-hand side anomalies are unable to be identified in the result by the linearized EIT. Also for the extended target, the proposed method shows better performance in reconstructing the conductivity value as well as the anomaly shape compared to that of the linearized EIT. Next, we show how the performance is degraded as the measurement points geometry changes. It is worth to remark that for $m=32$, which is a practical number of measurement points, the results are similar to those with $m=100$. Even when $m=16$, the obtained image shows the anomalies more clearly than the linearized EIT with $m=100$. {However, when $m=16p$, the results are distorted especially for the case of extended target.} In the perspective of the maximum recoverable targets related with $\mbox{spark}(A)$, we summarized the upper bound of the rightmost term in \eqref{eq:y=Ax0MMVuniq} for various number of measurement points in Table.\;\ref{tbl:UpperBound} with the value of $\|X\|_0$ for the sparse targets. Here, $\|X\|_0$ indicates the twice of the real number of non-zero positions of anomalies due to the pairwise joint sparsity structure in the EIT problem as described in Fig.\;\ref{fig:CS}(b). Note that $\|X\|_0$ of sparse target A is less than the upper bound for $m=100$ and larger than the one for $m=32$. However, it still shows a good reconstruction result when $m=32$, which implies that the theoretical upper bound using $\mbox{spark}(A)$ is rather conservative in practice where the exact recovery is impossible.
\begin{table}[!htbp]
\vspace*{0.2cm}
\centerline{
\renewcommand{\arraystretch}{1.3}
\begin{tabular}{|c|c|c|c|c|} \hline
 & $\|X\|_0$ & m=100 & m=32 & m=16\\
\hline\hline
 Sparse target A  & 48 & \multirow{2}{*}{51} & \multirow{2}{*}{17} & \multirow{2}{*}{9}\\
 Sparse target B  & 56 &                     &                     &                   \\
\noalign{\hrule height 0.5pt}
\end{tabular}
}\caption {The value of $\|X\|_0$ when $\Omega$ is discretized to have a grid size $0.5$, and the upper bound of the maximum recoverable positions for various number of measurement points.}
\label{tbl:UpperBound}
\end{table}

\begin{figure}[!htbp]
\centering\includegraphics[width=15cm]{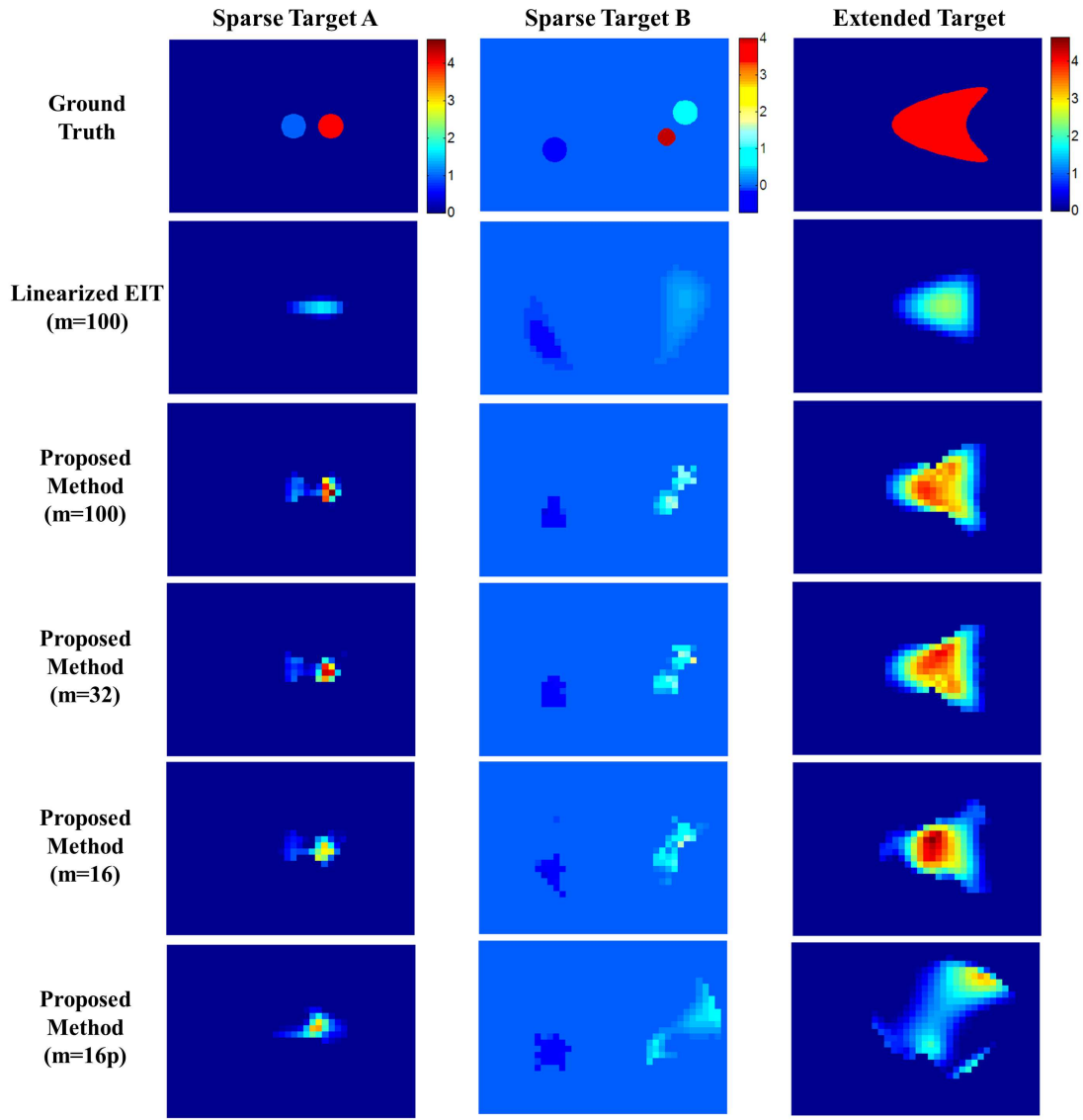}
\caption{Reconstruction results ($\sigma-1$) using the linearized EIT ($m=100$) and proposed method for various measurement points geometry with $M=2$.}
\label{fig:ResultAll}
\end{figure}

The relative error is calculated based on \eqref{eq:relError} for sparse and extended targets and is summarized in Table.\;\ref{tbl:ResultsError}. The value in the table denotes the average relative error from $20$ trials and the value in the parenthesis is their standard deviation. As we can see, the proposed method shows lower error compared to the linearized EIT for both sparse and extended targets.

\begin{table}[!htbp]
\vspace*{0.2cm}
\centerline{
\renewcommand{\arraystretch}{1.3}
\begin{tabular}{|c|c|ccc|}
\noalign{\hrule height 0.5pt}
         &         & Sparse target A   & Sparse target B   & Extended target \\\hline \hline
Linearized EIT & m=100 & 0.6368 ($\pm$ 0.0087) & 0.8909 ($\pm$ 0.0008) & 0.5096 ($\pm$ 0.0046) \\
\hline
\multirow{4}{*}{Proposed Method}& m=100 & 0.3697 ($\pm$ 0.1166) & 0.6296 ($\pm$ 0.0561) & 0.2227 ($\pm$ 0.0123) \\
 & m=32 & 0.4304 ($\pm$ 0.1260) & 0.7351 ($\pm$ 0.1141) & 0.2644 ($\pm$ 0.0286) \\
 & m=16 & 0.4705 ($\pm$ 0.1123) & 0.8577 ($\pm$ 0.1019) & 0.3592 ($\pm$ 0.0829) \\
 & m=16p & 0.6507 ($\pm$ 0.0344) & 0.9685 ($\pm$ 0.0471) & 0.6796 ($\pm$ 0.0391) \\
\noalign{\hrule height 0.5pt}
\end{tabular}
}\caption{Relative error for the reconstruction results with $M=2$.}\label{tbl:ResultsError}
\end{table}

The average reconstruction time (in {\rm [sec]}) of various methods for the sparse and extended target simulations are summarized in Table\;\ref{tbl:RunTime} (using a PC with CPU : core i7 sandy bridge). As is illustrated in Fig.\;\ref{fig:Compare}, the dimension reduction of the proposed method makes the run time of calculating the conductivity of the anomalies (using the C-SALSA algorithm) much faster than that of using the matrices obtained using the {linearized EIT}. Note that the proposed method has additional steps of finding non-zero support using M-SBL algorithm and estimating the unknown internal potential $\hat{u}$. However, the run time of M-SBL algorithm depends on the number of measurement points whose dimension is much smaller than that of the area of our interest, so this step is quite fast as described in Table\;\ref{tbl:RunTime}, and estimating $\hat{u}$ requires only simple matrix multiplication as described in \eqref{eq:ukxUkxInt2}. Therefore, the total run time of the proposed method is faster than that of the linearized EIT. The difference between total run time and C-SALSA (and M-SBL) is mostly dedicated for generating the sensing matrix.

\begin{table}[!htbp]
\vspace*{0.2cm}
\centerline{
\renewcommand{\arraystretch}{1.3}
\begin{tabular}{|c|cc|c|c|c|c|} \hline
 & \multicolumn{3}{c|}{Proposed method} & \multicolumn{2}{c|}{Linearized EIT}\\
\cline{2-6}
& \footnotesize{M-SBL} & \footnotesize{C-SALSA} & Total & \footnotesize{C-SALSA} & Total\\
\hline\hline
 Sparse target A  & 0.145 & 0.103 & 0.705 & 0.830 & 1.290\\
 Sparse target B  & 0.153 & 0.092 & 0.688 & 1.777 & 2.227\\
 Extended target  & 0.135 & 0.117 & 0.724 & 1.349 & 1.801\\
\noalign{\hrule height 0.5pt}
\end{tabular}
}\caption {Average run time [{\rm sec}] for anomalies reconstruction with $m=100$ and $M=2$.}
\label{tbl:RunTime}
\end{table}
Finally, we compare the proposed method to the MUSIC algorithm for localization of the anomalies. For this test, we take the absolute values of the reconstructed $(\sigma-1)$ using the proposed method and normalize it by scaling between $0$ and $1$. Accordingly, the calculated MUSIC spectrum based on \eqref{eq:MUSIC} is normalized. Fig.\;\ref{fig:ResultMUSIC} shows that MUSIC fails to localize the left anomaly for the case of sparse target A and the centered one for sparse target B when $M=4$. However, the proposed method can distinguish each anomalies even when $M=2$.

\begin{figure}[!htbp]
\centering\includegraphics[width=15cm]{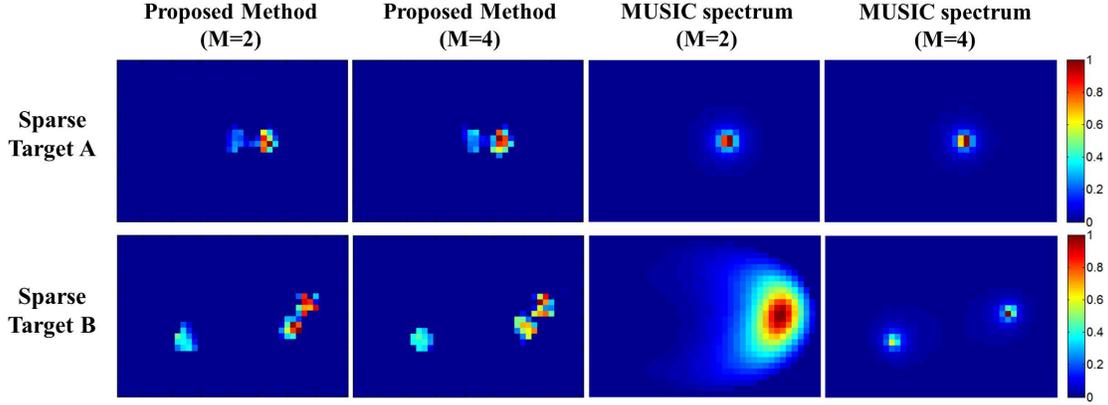}
\caption{The MUSIC spectrum and the normalized $|\sigma-1|$ obtained by the proposed method ($m=100$).}
\label{fig:ResultMUSIC}
\end{figure}

\section{Discussion and Conclusion}
\label{sec:Conclusion}

If the background conductivity is not a constant but a smooth function, say $\sigma_0$, then the background solution $U_k$ is the solution to \eqnref{eqn:u_k} with $\sigma_0$ in the place of $\sigma$. Using the variational formulation of $u_k$ and $U_k$, one obtains
$$
\int_\GO \sigma_0\nabla(u_k-U_k)\cdot\nabla v\;dy=\int_{\cup_{p=1}^N D_p}(\sigma_0-\sigma)\nabla u_k\cdot\nabla v\;dy\quad\mbox{for }v\in H^1(\GO).
$$
By substituting the Neumann function corresponding to $\sigma_0$, which we still denote by $N(x,y)$, to $v(y)$, we have the following equation:
$$u_k(x)-U_k(x) = \int_{\cup_{p=1}^N D_p} (\Gs_0(y)-\Gs(y))\nabla_y N(x,y) \cdot \nabla u_k(y)~ dy, \quad x \in \GO.$$
We emphasize that this formula is the same as \eqnref{ukxUkx} when the background conductivity $\sigma_0$ is constant. Hence the method proposed in this paper can be applied to reconstruct small anomalies embedded in the background domain with the continuous conductivity if one has the background solution and the Neumann function for the background domain.

Let us further consider the EIT system using $M$ number of electrodes $\mathcal{E}_k$ for $k=1,\dots,M$. Let $u_k$ be the potential subject to the $k$-th pairwise injection current $I$ between adjacent electrodes $(\mathcal{E}_k, \mathcal{E}_{k+1})$ with $\mathcal{E}_{M+1}=\mathcal{E}_1$. Ignoring the effects of the unknown contact impedances between electrodes and $\partial\Omega$, $u_k$ approximately satisfies the Neumann boundary problem \eqnref{eqn:u_k} with $g_k=0$ on $\GO\setminus(\mathcal{E}_k\cup\mathcal{E}_{k+1})$ and 
$\int_{\mathcal{E}_k}g_k\;ds=I=-\int_{\mathcal{E}_{k+1}}g_k\;ds$ for all $k=1,\dots,M$. 
We measure the boundary voltage difference $V_{j,k}[\Gs]$ between two adjacent electrodes $(\mathcal{E}_j, \mathcal{E}_{j+1})$, which is 
$$V_{j,k}[\Gs]=\frac{1}{|\mathcal{E}_j|}\int_{\mathcal{E}_j}u_k\; ds- \frac{1}{|\mathcal{E}_{j+1}|}\int_{\mathcal{E}_{j+1}}u_k\; ds,$$ 
then the following equation can be derived (see Ch 7 in \cite{SWbook}):
\beq\label{eqn:electrode}
V_{j,k}[\Gs]-V_{j,k}[\Gs_0]=\int_{\cup_{p=1}^N D_p} (\sigma_0-\sigma)\nabla u_k\cdot \nabla U_j\; dy\quad\mbox{for }j,k=1,\dots,M.
\eeq
Here $U_j$ and $V_{j,k}[\Gs_0]$ are respectively the background potential and the voltage difference corresponding to the background conductivity $\sigma_0$. Since the support of the integral in \eqnref{eqn:electrode} is independent of $j,k$, one can apply the joint sparsity recovery to obtain $(\sigma_0-\sigma)\nabla u_k$. Unfortunately, \eqnref{eqn:electrode} is valid only on the boundary of $\GO$ and, thus, $\nabla u_k$ cannot be computed. But, one may apply C-SALSA algorithm with $\nabla U_k$ instead of $\nabla u_k$ to compute the conductivity difference $\sigma_0-\sigma$.

This paper proposed the reconstruction method that resolves the non-linearity of the inverse problem of electrical impedance tomography. It accurately reconstructs the anomalies without iterative procedure or linear approximation. The main idea of the proposed method comes from the joint sparsity in the compressed sensing theory. The non-linear inverse problem of EIT can be changed into the joint sparse recovery problem, and it enables us to estimate the unknown internal potential with the help of the recursive nature of the forward problem formulation. Finally, the electrical property of the anomalies can be calculated from the proposed linear problem. The simulation results showed that the proposed method outperforms over the linearized EIT method {and the MUSIC algorithm}. 
Restriction of the region of interest to the estimated positions of anomalies and the estimation of the unknown internal potential enable the proposed method to reconstruct the anomalies in more fast and accurate way compared to the linearized method. It can also discriminate sparse anomalies which are located with a distance comparable to the anomaly size while MUSIC fails due to the lack of number of independent measurements. On the other hand, it faces to calculate the Neumann function when the measurement points are available only on the part of the boundary, and the estimation of the internal potential is limited in the complete electrode model.


\end{document}